\documentclass[12pt]{amsart}
\usepackage{amssymb, amscd, amsmath, latexsym, 
}
\newtheorem{thm}{Theorem}[section]
\newtheorem{cor}[thm]{Corollary}
\newtheorem{defn}[thm]{Definition}

\newtheorem{prop}[thm]{Proposition}
\newtheorem{rem}[thm]{Remark}
\numberwithin{equation}{section}
\newtheorem{exl}[thm]{Example}

\numberwithin{equation}{section}

\def\proof{\textsc{Proof:\ }}
\def\endproof{$\Box$\medskip}

\def\S+{\mathbf{SO(n, 1)^+}}
\def\O+{\mathbf{O(n, 1)}}

\def\s+{\mathcal{S}}

\def\proof{\textsc{Proof:\ }}
\def\endproof{$\Box$\medskip}

\begin{document}
\title{Twisted Automorphisms of Right Loops}
\author{R. Lal}
\address{IIIT Allahabad (INDIA)\ 211002\\
\textnormal{ mathrjl@gmail.com}}

\author{A. C. \ Yadav}
\address{Department of Mathematics\\
 M G Kashi Vidyapith, Varanasi (INDIA)\\
\textnormal{akhileshyadav538@gmail.com}}

\date{  August 5, 2013}
\begin{abstract}
In this paper we make an attempt to study right loops $(S, o)$ in which, for each $y\in S$, the map $\sigma_y$ from the  inner mapping group $G_S$ of $(S, o)$ to itself given by $\sigma_y (h)(x) o\ h(y)= h(xoy)$, $x\in S, h\in G_S$ is a homomorphism. The  concept of twisted automorphisms of a right loop and also the concept of  twisted right gyrogroup appears naturally and it turns out that the study is almost equivalent to the study of  twisted automorphisms and a twisted right gyrogroup. A representation theorem for twisted right gyrogroup is established. We also study relationship between twisted gyrotransversals and twisted  subgroups (a concept which arose as a tool to study computational complexity involving class NP). 
 \end{abstract}
\maketitle
{\bf Keywords}: Gyrotransversals, Right gyrogroups, Inner mapping groups, Twisted subgroups, CSP.

{\bf MSC}: 20N05, 08A35, 20D45.

\section{Introduction}

It has always been practice to study right loops through its  right inner mappings and  inner mapping groups (also called group torsions~\cite{lal}). However, for each right loop $(S, o)$,  there is another important family $\{\sigma_y : G_S \rightarrow G_S\ |\ y \in S\}$ of maps from the inner mapping group $G_S$ of $(S, o)$ to itself given by 
\begin{eqnarray*}
h(xoy) & =  & \sigma_y(h)(x) o h(y),\  h \in G_S
\end{eqnarray*}
which has not been so explored. 
Of course, right loops for which each $\sigma_y$ is identity map (equivalently all members of $G_S$ are automorphisms) termed as $A_r -loops$ has been studied (\cite{bruck, fog:ung}). In this paper we initiate the  study of right loops for which the family $\{\sigma_y :\  y \in S\}$ has some prescribed properties. Our special stress here will be on right loops for which the maps $\sigma_y, y\in S$ described above are automorphisms. In turn, we come across the concept of a twisted automorphism and twisted right gyrogroup. We have a representation theorem of twisted right gyrogroup as twisted gyrotransversal (Theorem~\ref{mainth}). 

The concept of near subgroups of a group was introduced by Feder and Vardi~\cite{fed:vardi} in an attempt to study CSP and which can be defined as twisted subgroups with some additional property.
  Indeed, twisted subgroups and  near subgroups in  a group of odd order are same, and also order dividing twisted subgroups are always near subgroups~\cite{fed}. The structure theory of near subgroups was ingeniously developed by M. Aschbacher~\cite{aschbacher}. Twisted subgroups, therefore, are important objects to study complexity of constraint satisfaction problems (CSP). 
Foguel and Ungar (\cite{fog:ung}, Theorem 3.8) have described equivalence of twisted subgroups and gyrotransversals under certain condition. Here in this paper we show [Theorem~\ref{twsgpgyro}]  that a similar equivalence holds between twisted gyrotransversals and twisted subgroups.

\section{Preliminaries}

Let $(S, o)$ be a right loop with identity $e$ and $y, z$ in $S$. 
 The map $f(y, z)$ from $S$ to $S$ given by the equation 
\begin{eqnarray}\label{bas1}
f(y,z)(x)o (yoz) & = & (xoy)oz,\ \qquad x\in S
\end{eqnarray}
belongs to $Sym\ S$ (the Symmetric group on $S$) and is called a right inner mapping of $(S, o)$.
Indeed $f(y, z)\in Sym (S\setminus\{e\})\subseteq Sym\ S$.
The subgroup $G_S$ of $Sym (S\setminus\{e\})\subseteq Sym\ S$ generated by $\{f(y, z)\ | \, y, z\in S\}$ is  called the inner mapping group (also called  the group torsion \cite{lal}) of $(S, o)$.

Further, let $h\in Sym (S\setminus\{e\})\subseteq Sym\ S$ and $y\in S$.  
Define $\sigma_y (h)\in Sym (S\setminus\{e\})\subseteq Sym\ S$ by the equation
\begin{eqnarray}\label{rlal}
h(x o y)  & = &  \sigma_y (h) (x) o h(y),  \  \qquad x\in S
\end{eqnarray} 

 \begin{prop}\label{rlal1}\cite{lal}
Let $(S, o)$ be a right loop with identity $e$. Then we have the following identities:\\
(i)\  $f(x, e)= f(e,x) = I_S$, for each $x\in S$.\\
(ii)\  $\sigma_e = I_S$, the identity map on $S$.\\
(iii)\   $\sigma_y (hk) =\sigma_y (h)\sigma_{ h(y)} (k)$, for all $y\in S$ and $h, k\in G_S$. In particular $\sigma_x (I_S) = I_S$.\\ 
(iv)\  $f(x, y) f(xoy, z) = \sigma_x (f(y,z)) f( f(y, z) (x), yoz)$, for all $x, y, z\in S$.\end{prop}
\begin{prop}\label{rlal11}
Let $(S, o)$ be a right loop with identity $e$. Suppose that $f (x', x)= I_S$ for all $x\in S$, $x'$ being the left inverse of $x$. Then,\\
(i)\ each element $x\in S$ has the unique inverse $x'$ in the sense that $x' ox= xox'= e$. In particular, it is a right loop with unique inverses.\\
(ii)\ $f(x, x')= I_S$. 
\end{prop}
\proof
(i) \begin{eqnarray*}
(xo x') o x & = &  f(x', x)(x) o (x' o x) \qquad {\rm (by\ equation~\ref{bas1}  )}\\
& = & xo e\\
& = & x\\
& = &  eox
\end{eqnarray*}
 By right cancellation law in right loop, $x o x'= e$.

(ii)\ By Proposition  \ref{rlal1}(iv),
\begin{eqnarray*}
f(x, x') f(x o x', x)& = & \sigma_x (f(x', x)) f( f(x', x)(x), x' ox)
\end{eqnarray*}
By Proposition~\ref{rlal11} (i) and Proposition~\ref{rlal1} (i), we get $f(x, x')= I_S$.
\section{Twisted Automorphisms}

\begin{defn}\label{thom} (Twisted automorphisms)\ 
Let $(S, o)$ be a right loop together with  unique inverse $x'$ for each $x$ ($x' o x= e= x o x'$). Then a bijective map $h:S\rightarrow S$ is called a twisted automorphism if 
\begin{eqnarray}\label{t1}
h(xoy) & =  & [h(x')]' o h(y)
\end{eqnarray} 
for all $x,y\in S$ with $y\neq e$.  \end{defn}

 For the sake of convenience  we shall also write $x\theta h$ for $h(x)$. Thus the equation (\ref{t1}) also reads as 
\begin{eqnarray}\label{t2}
(xoy)\theta h & =  & [(x')\theta h]' o y\theta h.
\end{eqnarray} 

As we shall be dealing with right transversals  and right actions, it is also convenient for us to adopt the convention $(gf) (x)= f(g(x))$ for composition of maps. Thus, $x\theta (gf)= (x\theta g)\theta f$ for all maps $f,g$ on $S$  and $x\in S$. 
\begin{rem}
(i)\  If $h$ is a twisted automorphism of $(S, o)$, then  $h(e)= h(x' o x)= (h(x''))' o h(x)= (h(x))' o h(x) = e$, $x\neq e$.\\
(ii)\  If $h$ is an automorphism of $(S, o)$, then $h(x) o h(x')=  h(xo x')=h(e)= e$. By uniqueness of inverse,  $h(x) = [h(x')]'$ and so  $h$ is a  twisted automorphism of $(S, o)$.\\
(iii)\  A twisted automorphism of $(S, o)$ need not be an automorphism of $(S, o)$ (c.f. Examples\ \ref{ex1}, \ref{4ag}).
 \end{rem}
 
\begin{prop}
Let $(S, o)$ be a right loop together with  unique inverse $x'$ for each $x$.  Then the set $TAut\ (S, o)$ of all twisted automorphisms of the right loop $(S, o)$ forms a group under the composition of maps. 
\end{prop}
\proof
Let $h, k\in TAut\ (S, o)$. Then, for $y\ne e$ 
\begin{eqnarray*}
(xoy)\theta (hk) & = & \left((xoy)\theta h\right)\theta k\\
& = & \left[((x'\theta h)') o (y\theta h)\right] \theta k\\
& = & \left[((((x'\theta h)')'\theta k)') o ((y\theta h)\theta k)\right] \qquad  \ {\rm as}  \, y\theta h\ne e\\
& = & \left[(x'\theta h)\theta k\right]' o (y\theta (hk))\\
& = & (x'\theta (hk))' o (y\theta (hk))
\end{eqnarray*}
 Thus,  $hk\in TAut\ (S, o)$.  Clearly $I_S\in TAut\ (S, o)$. Let $h\in TAut\ (S, o)$ and $u, v\in S$ with $v\ne e$. Let  $x, y\in S$ with $y\ne e$ such that $x\theta h = u, \ y\theta h = v$. Then
\begin{eqnarray*}
(x'oy)\theta h& =  & [(x'')\theta h]' o y\theta h\\
&  = & (x\theta h)' o y\theta h\\
& = & u' o v
\end{eqnarray*}
Thus, 
\begin{eqnarray*}
(u' o v)\theta h^{-1}  & = & x'o y\\
& = & (u\theta h^{-1})' o v\theta h^{-1}
\end{eqnarray*} 
Replacing $u$ by $u'$, we have 
$$(u o v)\theta h^{-1} = (u'\theta h^{-1})' o v\theta h^{-1}.$$
 This shows that  $h^{-1}$ is also a  twisted automorphism and so  $TAut\ (S, o)$ is a group.
\endproof

\begin{prop}
Every twisted automorphism of a group is an automorphism. \end{prop}
\proof
Let $G$ be a group.  Let $h\in TAut (G)$.  To prove that $h\in Aut\ G$, it is sufficient  to prove that $(h(x'))'= h(x)$, or equivalently $h(x')= (h(x))'$.\\
{\bf Case I:} If $\left|G\right|= 2 $, then $h= I_G$ and there is nothing to do.\\
{\bf  Case II:} If $\left|G\right|= 3 $, then for $x\ne e$, $x^2 = x' \ne e$. Further, 
\begin{eqnarray*}
h(x) & = & h(x^2. x^2)\\
& = & [h((x^2)')]'. h(x^2)\\
& = & h(x)'  h(x') \ as \ (x^2)' = x
\end{eqnarray*}
Thus, $h(x)^2 =  h(x)h(x)'h(x')= h(x')$. This shows that $h(x)' =h(x')$.\\
{\bf Case  III:} Suppose that $\left|G\right|\ge 4$. Let $x,y,z\in G\setminus \{e\}$ with $yz\ne e$. Then, since $h$ is a twisted automorphism,  
\begin{eqnarray*}
h[(xy)z] & = & [h((xy)')]' h(z)\\
& = & [h(y' x')]'h(z)\\
& = & [h(y)' h(x')]' h(z)
\end{eqnarray*}
and 
\begin{eqnarray*}
h[x(yz)] & = & h(x')' h(yz)\\ 
& = & h(x')' [h(y')' h(z)] \\
& = & [h(y') h(x')]' h(z)
\end{eqnarray*}
Since  $h[(xy)z] = h[x(yz)]$, $[h(y)' h(x')]  = [h(y') h(x')]$ and so $h(y') = h(y)'$ for all $y\ne e$.  
\endproof

Again,  we have the following:
\begin{prop}
Let $(S, o)$ be  a right loop together with unique inverses and which  satisfies the following identities:\\
(i)\  $(ab)' = a' b'$ (automorphic inverse property (AIP)), \\
(ii)\  $(aa)b  = a(ab)$ (left alternative)\\
 for all $a, b\in S$. Then every twisted automorphism of $(S, o)$ is an automorphism. 
\end{prop}

\proof 
Assume  the conditions (i) and (ii). Let $h$ be a twisted automorphism. Then, since every right loop containing two elements is a group, it is sufficient to assume that $S$ contains more than two elements. Let $a\in S\setminus\{e\}$. Then $\exists \ b\neq e$ such that $ab\neq e$ for, otherwise, each element of $S\setminus \{e\}$ will be a inverse of $a$, a contradiction to the assumption that $(S, o)$ is with unique inverses. Thus
\begin{eqnarray*}
h[a(ab)] & = & h(a')'. h(ab)\\
& = & h(a')'.[ h(a')' h(b)]\\
& = & [h(a')'. h(a')'] h(b) \ \qquad \qquad{\rm( by \ (ii))}\\
& & \\
h[(aa)c] & = & h((a^2)')'. h(b)\\
& = & [h(a'. a') ]'. h(b)\ \qquad \qquad {\rm ( by \ (i))}\\
& = & [h(a'')'. h(a')]' h(b)\\
& = & [h(a)'. h(a')]' h(b)\\
& = & [h(a). h(a')']. h(b) \qquad \qquad {\rm (by \ (i))}
\end{eqnarray*}
Since  $h[a(ab)]= h[(aa)b]$,  $[h(a')'. h(a')']= [h(a). h(a')']$ and so $h(a')'= h(a)$ for all $a\in S$.
\endproof

\begin{cor}
Every twisted automorphism of a K-loop  is an automorphism. $\Box$
\end{cor}

\section{Twisted Right Gyrogroups And Twisted Gyrotransversals}
In this section, we consider right loops $(S, o)$ in which 
for each $y\in S$, the map $\sigma_y$ from $G_S$ to $G_S$ given by the equation~(\ref{rlal}) is an automorphism and discuss its equivalence with the requirement that all right inner mappings of $(S, o)$ are twisted automorphisms.
\begin{prop} \label{fundam}
Let $(S, o)$ be  a  right loop.  Suppose that  $G_S$ acts transitively on $S\setminus \{e\}$ and the map  $\sigma_y :G_S\rightarrow G_S$  defined by
\begin{eqnarray}\label{eta1}
x\theta \sigma_y(h) o y\theta h & =  & (xoy)\theta h
\end{eqnarray} 
is a homomorphism for some $y\in S\setminus \{e\}$. Then $\sigma_x =\sigma_y$ for each $x\in S\setminus \{e\}$. If in addition $f(x', x)= I_S$ ($x'$ being the left inverse of $x\in S$), then $\sigma_y$ is an involutory automorphism of $G_S$. Further, in this case $G_S\subseteq TAut\ (S, o)$.
\end{prop}
\proof 
Let $x,y \in S\setminus \{e\}$ with $\sigma_y$  a homomorphism. By transitivity of the action of $G_S$ on $S\setminus \{e\}$, there exists $h\in G_S$ such that $y\theta h = x$. Let $k\in G_S$. Then 
\begin{eqnarray*} \label{dtwist2}
\sigma_y (h) \sigma_y (k) & = & \sigma_y (hk) \\
& = &  \sigma_y (h) \sigma_{y\theta h} (k)\qquad ({\rm Proposition }~\ref{rlal1} (iii))\\
& = & \sigma_y (h) \sigma_x (k)
\end{eqnarray*} 
This gives that $\sigma_y (k)= \sigma_x (k) $ for each $k\in G_S$. In turn,  $\sigma_x = \sigma_y$ for each $x\in S\setminus \{e\}$. Denote $\sigma_y $ by $ \eta$. Then $\sigma_x = \eta$ for each $x\in S\setminus \{e\}$.  
Now, the equation (\ref{eta1}) reduces to
\begin{eqnarray}\label{eta2}
x\theta \eta (h) o y\theta h & =  & (xoy)\theta h
\end{eqnarray} 
 and the  identity $(iv)$ of Proposition~\ref{rlal1} becomes\\
\begin{eqnarray}\label{pp1}
f(x, y) f(xoy, z) & = & \eta (f(y,z)) f(x\theta f(y, z), yoz)
\end{eqnarray} 
Assume that $f(x', x)= I_S$.  By Proposition~\ref{rlal11},   $(S, o)$ is a right loop with unique inverses. Using equation~(\ref{pp1}), we have 
$$f(y', y) f(y'oy, z) = \eta (f(y, z)) f(y'\theta f(y, z), yoz). $$
This gives 
\begin{eqnarray} \label{involut}
\eta (f(y,z)) & = {f(y'\theta f(y, z), yoz)}^{-1}\,&  \forall  y, z\in S\setminus\{e\}.
\end{eqnarray}
and so 
\begin{eqnarray}\label{inv1}
\eta^2 (f(y, z)) & = & \eta \left(f( y'\theta f(y, z),  yoz )^{-1}\right)\nonumber\\ 
& = & [\eta \left(f(y'\theta f(y, z), yoz)\right)]^{-1} \nonumber\\
& = & \{\eta \left(f(Y, Z)\right)\}^{-1} \, \textrm{where}\ Y =  y'\theta f(y, z) , \, Z = yoz\nonumber\\
 & = & f(Y'\theta f(Y, Z), Yo Z)\qquad ( {\rm by\ equation\  (\ref{involut}}))
\end{eqnarray}
   Now 
   \begin{eqnarray*}\label{inv2}
   YoZ & = &  y'\theta f(y, z)o (yoz) \\
   & = &  (y'oy)oz\\
   & = & eoz\\
   & =&  z
   \end{eqnarray*}
   and also 
   \begin{eqnarray*}
   e   & = & (yoy')\theta f(y, z)\qquad \qquad (y'\ne e) \\
   & = & y\theta\eta( f(y, z))o y'\theta f(y, z)\qquad ( {\rm by\ equation\  (\ref{eta2}}))\\
   & = & y\theta{f(y'\theta f(y, z), yoz)}^{-1} o Y
\end{eqnarray*}
Hence 
  $$Y' =  y\theta{f(y'\theta f(y, z), yoz)}^{-1}= y\theta f(Y, Z)^{-1};$$  
   and so \begin{eqnarray}\label{inv3}
	Y'\theta f(Y, Z) & =& y.
	\end{eqnarray}
 
 Thus, from equation (\ref{inv1})   $\eta^2 (f(y, z)) = f(y,z)$ for all  $y,z\ne e$. But $G_S = \left\langle f(y, z)|\, y,z\in S\right\rangle$ and so $\eta^2 = I_{G_S}$. This shows that  $\eta$ is an involutory automorphism of $G_S$.
Next, let $x\in S\setminus \{e\}$. Since $x'$ is also a right inverse of $x$, by equation (\ref{eta2})
\begin{eqnarray*}
e & = & (xox')\theta f(u, v) \\
& = &  x\theta \eta [f(u, v)] o x'\theta f(u, v) 
\end{eqnarray*}
 Thus $ x\theta \eta [f(u, v)] = [x'\theta f(u, v)]' $ for all $x\in S$ and so
\begin{eqnarray*}
(xoy)\theta f(u, v) & = & x\theta \eta [f(u, v)] o y\theta f(u, v)\\
& = & [x'\theta f(u, v)]' o y\theta f(u, v).
\end{eqnarray*}
This shows that $f(u, v)$ is a twisted automorphism of $(S, o)$ for all $u, v $ in  $S$. Since $TAut\ (S, o)$ is a group under composition of maps, it follows that $G_S$ is a subgroup of $TAut\ (S, o)$.
\endproof

\begin{cor}\label{20jul}
Let $(S, o)$ be a loop such that $\sigma_y : G_S\rightarrow G_S$ defined by equation~(\ref{rlal}) is an automorphism for some $y\in S\setminus\{e\}$. Then $\sigma_x = \sigma_y= \eta$ is an involutory automorphism of $G_S$ and all members of $G_S$ are twisted automorphisms.
\end{cor}

\proof
Since $(S, o)$ is a loop, $G_S$ acts sharply transitively on $S\setminus \{e\}$(\cite{kiec}, Theorem 2.11(2)). The result follows from Proposition ~\ref{fundam}.
\endproof

\begin{rem}\label{4agr}
 Note that even if $(S,o)$ is not a loop, $G_S$ may act transitively on $S\setminus\{e\}$ (c.f. Example~\ref{ex1}).
\end{rem}

The above Proposition prompts us to have  the following:

\begin{defn}(Twisted Right Gyrogroups)\label{trgyrogp}
A right loop $(S, o)$  is said to be a twisted right gyrogroup if \\
(i)\  $f(y', y) = I_S$ for all $y\in S$ (Consequently (S, o) is a right loop with unique inverses),\\
(ii)\  the right inner mappings  $f(y, z): S\rightarrow S$  are twisted automorphisms for all $y,z$ in $S$.  
\end{defn}
The following Proposition is a partial converse of the Proposition~\ref{fundam}. 
\begin{prop}\label{20julp}
Let $(S, o)$ be a twisted right gyrogroup. Then it determines an involutory automorphism $\eta : G_S\rightarrow G_S$ such that 
\begin{eqnarray}\label{X1}
(xoy)\theta h & = & x\theta \eta(h) o y\theta h
\end{eqnarray}
 \ for\  $x,\ y\in S$ and $h\in G_S$. In particular all $\sigma_y\ ( = \eta)$ are automorphisms.
\end{prop}
\proof
Let $x\in S$ and $h\in Sym\ (S\setminus \{e\})\subseteq Sym\ S$. Then the map $x\rightarrow (x'\theta h)'$ is bijective. This determines a map $\eta : Sym\ (S\setminus \{e\})\rightarrow Sym\ (S\setminus \{e\})$ given by 
\begin{eqnarray}\label{definvo}
x\theta \eta(h) & = & (x'\theta h)'
\end{eqnarray}
Let $h, k\in Sym\ (S\setminus \{e\})$ and $x\in S$. Then 
\begin{eqnarray*}
x\theta \eta (h.k) & = & [x'\theta (h.k)]'\\
& =  &  ((x'\theta h)\theta k)'\\
& = & (((x'\theta h)')'\theta k)'\\
& = & ((x'\theta h)')\theta \eta (k)\\
& = & (x\theta \eta (h))\theta \eta (k)\\
& = & x\theta (\eta(h). \eta (k))
\end{eqnarray*}
This shows that $\eta (hk)= \eta(h) \eta (k)$ and so $\eta$ is an endomorphism of $Sym\ (S\setminus \{e\})\subseteq Sym\ S$.

Next, let $x\in S$ and $h\in Sym\ (S\setminus \{e\})$. Then
\begin{eqnarray*}\label{Aut}
x\theta \eta^2(h)& = & x\theta\eta(\eta(h))\\
& = & (x'\theta\eta(h) )'\\
& = &  [((x'')\theta h)']'\\
& = & x\theta h\qquad \qquad\qquad ({\rm by\ uniqueness \ of\  inverses\  in }\ S)
\end{eqnarray*}
This shows that  $\eta^2 = I_{Sym\ (S\setminus\{e\})}$ and so $\eta$ is an involutory automorphism of $Sym\ (S\setminus \{e\})$.

Let $x, y, z\in S\setminus\{e\}$. Then
\begin{eqnarray*}
xoz & = & x o ((y' o y)o z)\\
& = & xo (YoZ), \ {\rm{where}\ Y= y'\theta f(y, z); Z= yoz}\\
& = & \left(x\theta f(Y, Z)^{-1} o Y\right) o (yoz)\\
&  = & \left(\left(x\theta f(Y, Z)^{-1} o Y\right)\theta f(y, z)^{-1} o y\right)oz\\
& = &  \left(\left(x\theta f(Y, Z)^{-1}\eta [f(y, z)^{-1}] o Y\theta f(y, z)^{-1}\right) o y\right)oz\\
& = &\left( \left(x\theta f(Y, Z)^{-1}\eta [f(y, z)^{-1}]f(Y\theta f(y, z)^{-1}, y)\right) o \left (Y\theta f(y, z)^{-1} o y\right)\right)oz
\end{eqnarray*}
But  $Y\theta f(y, z)^{-1}= y'$. Hence, $Y\theta f(y, z)^{-1} o y= y'oy = e$ and  $f(Y\theta f(y, z)^{-1}, y)= f(y', y)= I_S$. Using these, we have  
\begin{eqnarray*}
xoz & = &x\theta f(Y, Z)^{-1}\eta [f(y, z)^{-1}] o z
\end{eqnarray*}
By right cancellation in $S$, we have  $x\theta f(Y, Z)^{-1}\eta [f(y, z)^{-1}] = x $ for all $x$. This shows that
\begin{eqnarray}
\eta [f(y, z)] & = & f(Y, Z)^{-1} \end{eqnarray}
for all $y,z\in S$.  Thus, 
\begin{eqnarray}
\eta (G_S)&\subseteq& G_S
\end{eqnarray}
Now, let $h\in G_S$. Putting the value of $\eta (h),\ h\in G_S$  from equation (\ref{definvo}) in equation (\ref{X1}), we have 
\begin{eqnarray}\label{X2}
(xoy)\theta h & = & (x'\theta h)' o y\theta h.
\end{eqnarray}
This shows that $h\in TAut\ (S, o)$ for all $h\in G_S$. Further, $\sigma_y(= \eta)$ are automorphisms for each $y\in S$.
\endproof

\begin{defn}\label{tgyrotrans} (Twisted Gyrotransversals)
A right transversal $S$ to a subgroup $H$ of a group $G$ is said to be a twisted gyrotransversal if it satisfies the following:
\begin{enumerate}
\item $x^{-1}\in S$, for each $x\in S$,

\item there is an involution $\eta$ in $Aut\ H$ such that $\eta (h)^{-1} x h\in S$, for all $x\in S\setminus\{e\}$ and $h\in H$.
\end{enumerate} \end{defn}

\begin{prop}\label{lem1}
Let $(S, o)$ be a twisted right gyrogroup. Then  $S$ appears as a twisted gyrotransversal  to $G_S$ in a group $G$ such that the corresponding induced structure on $S$ is the given twisted right gyrogroup structure. 
\end{prop}
\proof
Since $S$ is a twisted right gyrogroup, by Proposition~\ref{20julp}, we have an involutory automorphism of $G_S$ such that $\sigma_x = \eta \  \forall x\neq e$.
Let us identify the order pair $(a, x)$ by $ax$ and the Cartesian product $G_S\times S$ by $G_SS$. Then $G_S S$ is a group with respect to the multiplication $\cdot$ given by 
\begin{eqnarray}
ax\cdot by & = & a\eta(b) f(x\theta b, y) ((x\theta b) o y),\ for \ x\neq e
\end{eqnarray} 
and $ae\cdot by= ab y$ for all $a, b$ in $G_S$, $x, y$ in $S$.  The identity of this group is $I_S e$.  The inverse of $ax$ is given by $$(ax)^{-1}= (f(x', x))^{-1} \eta (a^{-1}) (x'\theta a^{-1}), \ {\rm for} \ x\neq e$$ and $(ae)^{-1} = a^{-1} e$.
We identify $G_S$ with the subgroup $G_S \{e\}$ through the embedding $a\leadsto ae$ and $S$ with the subset $I_S  S$ of $G_S S$ through the embedding $x\leadsto I_S x$. Then it follows that $S$ is a right transversal to the subgroup $G_S$ of  $G_S S$. 
Let $x\in S$ and $h$ in $G_S$. Then 
\begin{eqnarray*}
x^{-1} & = & f(x', x)^{-1} x'= I_S x'= x'\in S
\end{eqnarray*} 
and 
$$x\cdot h = \eta (h) x\theta h,\ x\neq e .$$
 Thus, 
\begin{eqnarray*}
(\eta(h))^{-1}\cdot x\cdot h & = & x\theta h\in S,\ \forall\ x\ne e 
\end{eqnarray*}
and so  $S$ is a twisted gyrotransversal to the subgroup $G_S$ of the  group $G_S S$. Clearly
\begin{eqnarray*}
 x \cdot y & =& f(x, y)\cdot xoy
\end{eqnarray*}
  and so $S$ with the induced operation  is the given twisted right gyrogroup.
\endproof

Conversely,  we have the following
\begin{prop}\label{lem2}
Let $S$ be a twisted gyrotransversal to a subgroup $H$ in a group $G$. Then $S$ with the induced operation  is a twisted right gyrogroup. 
\end{prop}
\proof
Since $S$ is a right transversal to $H$ in $G$, we have a map $g:S\times S\rightarrow H$ and a binary operation $o$ on $S$ such that 
\begin{eqnarray}\label{31jul1}
xy & =&  g(x, y)\ xoy,\ x, y\in S
\end{eqnarray}

Clearly $(S, o)$ is a right loop. Further, since $x^{-1}\in S$ for every $x\in S$, it follows that the inverse $x'$ of $x$ in $(S, o)$ is given by $x' = x^{-1}$. Also  $g(x', x) =1= g(x, x')$. Again, we have an involutory automorphism $\eta$ of $H$ such that $\eta (h)^{-1} yh\in S$ for all $y\in S\setminus \{e\}$ and $h\in H$. Thus,
\begin{eqnarray}\label{31jul2}
yh= \eta(h)\ y\theta h,\ y\in S\setminus\{e\}, h\in H
\end{eqnarray} 
where $\theta$ is the right action of $H$ on $S$ given by $y\theta h= \eta(h)^{-1} yh, y\ne e$ and of course $e\theta h = e$ for every $h\in H$. Now, using equations (\ref{31jul1}), (\ref{31jul2}) and the fact that $\eta^2 = I_H$, we have for $y\ne e$, 
\begin{eqnarray*}
hg(x\theta \eta (h), y\theta h) \left (x\theta \eta (h) \ o \  y\theta h\right) & = & \eta (\eta (h)) \left[x\theta \eta (h)\cdot  y\theta h\right]\\
& = & [\eta (\eta (h))\  x\theta \eta (h)] \  (y\theta h)\\
& = & [x\  \eta (h)]\  y\theta h\\
& = & x (y h) \\
& = & (xy)h\\
& = & g(x, y) [(xoy)h] \\
& = & g(x, y)\eta (h) [(xoy)\theta h] 
\end{eqnarray*} 
This gives 
\begin{eqnarray}
hg(x\theta \eta (h), y\theta h)  & = & g(x, y)\eta (h), y\ne e\\
x\theta \eta (h)o y\theta h & = & (xoy)\theta h, y\ne e\label{tmain}
\end{eqnarray}

Next, if $x\ne e$, taking $y= x'$ in the  equation ~(\ref{tmain}), we get  
\begin{eqnarray*}
x\theta \eta (h) o x'\theta h & = &  (x o x')\theta h\\
& = & e\theta h\\
& = & e
\end{eqnarray*}
Hence 
\begin{eqnarray}
x\theta \eta (h) = [x'\theta h]' \ for\  x\ne e\label{tmain1}
\end{eqnarray}
Already $e\theta \eta(h) = e= (e'\theta h)'$.
Putting the value of $x\theta \eta(h)$ in the equation ~(\ref{tmain}), we get
$$(xoy)\theta h  =  [x'\theta h]'o y\theta h,\ for \ y\ne e$$
This shows that $H$ acts on $S$ through twisted automorphisms.
Clearly the permutation map on $S$ induced by $g(x, y)$ is the right inner mapping $f(x, y)$ of $(S, o)$. This ensures that $(S, o)$ is a twisted right gyrogroup. 
\endproof

Using Propositions \ref{lem1} and \ref{lem2}, we have 
\begin{thm}\label{mainth}(Representation theorem)
A right loop $(S, o)$ is a twisted right gyrogroup if and only if $S$ appears as a twisted gyrotransversal to a subgroup $H$ in a group $G$ inducing the given right loop structure.$\qquad \Box$
\end{thm}

\begin{prop}
Let $(S, o)$  be a finite twisted right gyrogroup. Suppose that  $G_S\cap Aut\ (S, o)= \{I_S\}$, then $G_S$ is an abelian group of odd order and in this case
\begin{eqnarray}\label{2ag}
(x o y)\theta h & = &x\theta h^{-1}\ o\ y\theta h,
\end{eqnarray} 
 for $h\in G_S$.
\end{prop}
\proof
By Proposition~\ref{20julp}, there is an involutory automorphism $\eta$ of $G_S$ given by   
\begin{eqnarray}\label{2ag1}
(xo y)\theta h & = & x\theta \eta(h) o y\theta h
\end{eqnarray}
for $x, y\in S$ and $h\in G_S$. From this it follows that 
\begin{eqnarray}\label{3ag}
h\in Aut\ (S, o)\Leftrightarrow \eta(h) = h
\end{eqnarray}
Suppose that $G_S$ is of even order. Since $\eta$ is an involutory automorphism, there exists $h\in G_S\setminus\{I_S\}$ such that $\eta(h)= h$. But then $h\in Aut\ (S, o)$.  This contradicts  $G_S\cap Aut\ (S, o)= \{I_S\}$. Thus $G_S$ is of odd order.

Next, consider the map $\phi : G_S\rightarrow G_S$ defined  by $\phi (h)= h^{-1}\eta (h)$. Suppose that $\phi (h_1) = \phi(h_2)$. Then $h_1^{-1}\eta(h_1)= h_2^{-1}\eta(h_2)$ and so $\eta(h_1h_2^{-1})=h_1 h_2^{-1}$. Since $G_S\cap Aut\ (S, o) = \{I_S\}$, $h_1 = h_2$. This shows that $\phi$ is injective and so also surjective. Let $h\in G_S$. 
Then $h= k^{-1} \eta (k)$ for some  $k\in G_S$ and so 
$$\eta (h)= \eta [k^{-1}\eta (k)]= \eta (k)^{-1} k= h^{-1}.$$ 
Thus the equation~(\ref{2ag1}) reduces to the equation~(\ref{2ag}). Since $\eta$ is an automorphism of $G_S$, $G_S$ is abelian. 
\endproof

\section{Some Examples}
In this section we give some examples of twisted automorphisms, twisted right gyrogroups and compute the twisted automorphism group $TAut\ (S, o)$, for some right loops.

\begin{exl}
 Every right gyrogroup is a twisted right gyrogroup.
\end{exl}

\begin{exl}
 Let  $S$ be a gyrotransversal  to the subgroup $H$ of  a group $G$.  Let $h$ be an involution of $H\setminus Z(H)$, $Z(H)$ denotes the center of $H$. Then the  map $\eta : H\rightarrow H$ defined by $\eta (k) = hkh$ is an involutory automorphism.
Consider the transversal $S_h = \{hx\, | \,x\in S\setminus \{e\} \}\cup \{e\}$. Then, for $x\ne e$,  $(hx)^{-1} = h. h^{-1}x^{-1}h\in S_h$ and $\eta (k)^{-1}hxk =hk^{-1} h.hxk= h.k^{-1}xk\in S_h$. Thus,  $S_h$ is a twisted gyrotransversal to the subgroup $H$ in a group $G$. Further, since $\eta(k)\ne k$, for some $k$, $S_h$ is not a gyrotransversal.  By Proposition~\ref{lem2}, the induced structure on $S_h$ is a  twisted right gyrogroup structure.
\end{exl}

\begin{exl}\label{ex1}
Consider the right loop $(S, o)$, where $S=\{1, 2, 3, 4, 5\}$ and the binary operation $o$ on $S$ is described by the following composition table:
\begin{displaymath}
\begin{array}{l|c|c|c|c|r|}
o& 1& 2& 3& 4& 5\\\hline
1& 1& 2& 3& 4& 5\\\hline
2& 2& 3& 1& 3& 3\\\hline
3& 3& 1& 2& 2& 2\\\hline
4& 4& 4& 4& 1& 4\\\hline
5& 5& 5& 5& 5& 1\\\hline
\end{array}
\end{displaymath}

From the composition table, it is clear that $(S, o)$ is a right loop with unique inverse. Observe that 
\begin{eqnarray*}
f(2, 2)= f(2, 3) = f(3, 2)= f(3, 3)= f(4, 4)= f(5, 5)=I_S,\\
f(2, 4)= f(4, 3)= (2 3 4), f(2, 5)= f(5,3)= (2 5 3),f(3, 4)= f(4, 2) = (2 4 3)\\
f(3, 5)= f(5, 2)= (2 5 3), f(4, 5)= f(5, 4)= (2 3)(4 5)
\end{eqnarray*}
One may easily verify that $(2 3 4), (2 3 5)$ and $(2 3) (4 5)$ are twisted automorphisms and so $G_S = \left\langle (2 3 4), (2 3)(4 5)\right\rangle\cong A_4 \subseteq TAut\ (S, o)$. Also $Aut\ (S,  o)= \{ I, (2 3), (4 5), (2 3)(4 5)\}$\ $\subseteq TAut\ (S, o)$  and $TAut\ (S, o)= \left\langle G_S, Aut\ S\right\rangle\cong S_4$. 
One may easily verify that the action of $G_S $ on $S\setminus\{1\}$ is transitive (though $(S, o)$ is not a loop) and the involutory automorphism $\eta$ (Proposition~\ref{fundam}), in this case, is the permutation of elements of $G_S\cong A_4$ represented by the cycle  decomposition    
$$\eta= ((234), (243)) ((235),(253))((245), (345))((254), (354))((24)(35), (25)(34)).$$
\end{exl}

\begin{exl}
Let $(S, o)$ be  a right gyrogroup. Let $\rho$ be an involutory automorphism of $Aut\ (S, o)$. Define a binary operation $o_\rho$ on $S$ by 
\begin{displaymath}
xo_\rho y  =  \left\{ 
\begin{array}{lcr}
\rho(x)oy= x\theta \rho \ o y &{\rm if} & y\neq e\\
x & {\rm if} & y=e
\end{array}
\right.
\end{displaymath}

 Then $(S,o_\rho)$ is a right loop together with unique inverse $\rho(x')= x'\theta \rho$ for each $x\in S$. Let us denote the left inverse of $x$ in $(S, o_\rho)$  by $x_\rho'$. Thus $x_\rho' = \rho(x')= x'\theta \rho$. 
Let $f^\rho (y,z)$ denote the right inner mapping of $(S, o_\rho)$, $y, z\in S$. For $y\ne z_\rho',\ z\ne e$, 
\begin{eqnarray*}
x\theta [f(y\theta \rho, z) \rho] o_\rho (y o_\rho z) & = & x\theta f(y\theta \rho, z)  o (y o_\rho z)\\
 & = & x\theta f(y\theta \rho, z)\  o\ (y \theta \rho\ o z)\\
  & = & (x\ o\ y\theta \rho) \ o\ z\\
	& = & (x\theta \rho\ o\ y)\theta\rho\  o\  z\qquad\ for \ \rho^2 = I_S\\
	& = & (x\ o_\rho\ y) o_\rho z
\end{eqnarray*}
It follows, therefore, that $f^\rho (y, z) = f(y\theta \rho, z)\rho$, $y, z\in S$ and $f^\rho (y, z)\in Aut\ (S, o)$. Thus, for $y\ne z_\rho'$  and  $u, v\in S$,   $v\neq e$
\begin{eqnarray*}
 (u o_\rho v)\theta f^\rho (y,z) & = &  \left(u\theta\rho\ o \ v\right)\theta f^\rho (y,z)\\
& = & (u\theta\rho)\theta f^\rho (y,z)\ o\  v\theta f^\rho (y,z)\qquad\ {\rm for}\ f^\rho (y, z)\in Aut\ (S, o)\\
& = & (u'\theta\rho)'\theta f^\rho (y,z)\ o\ v\theta f^\rho (y,z)\qquad {\rm for \ u''= u\ and\ \rho\in Aut\ (S, o)}\\
& = & (u'_\rho)'\theta f^\rho (y,z)\ o\ v\theta f^\rho (y,z)\\
& = & (u'_\rho\theta f^\rho (y,z))'\theta\rho\ o_\rho\  v\theta f^\rho (y,z)\qquad {\rm for}\ f^\rho (y, z)\in Aut\ (S, o)\\
& = & (u'_\rho\theta f^\rho (y,z))'_\rho\  o_\rho\  v\theta f^\rho (y,z)\\
\end{eqnarray*}
 This shows that $f^\rho (y,z)\in TAut\ (S, o_\rho)$.
Next,
 \begin{eqnarray*}
x\theta f^\rho (y_\rho', y)& = & x\theta f^\rho (y_\rho', y)\  o_\rho\  (y_\rho' o_\rho y) \\
& = & (x\ o_\rho\ y_\rho')\ o_\rho\ y\\
& = & (x\theta\rho\ o y'\theta \rho)\theta\rho \ o\ y\\
& = & (x o y') o y\qquad {\rm for} \ \rho\in Aut\ (S, o)\\
& = & x\theta f(y', y)o (y' o y) = x.
\end{eqnarray*}
This gives  that $f^\rho (y_\rho', y) = I_S$ and so $(S, o_\rho)$ is a twisted right gyrogroup.
\end{exl}

\begin{prop}\label{tat1}
 Let $(S, o)$ be  a right gyrogroup. Suppose that
 $Aut\ (S, o)$ contains an involution $\rho$. Then each $h\in Aut\ (S, o)$ is a twisted automorphism of $(S, o_\rho)$. Further, $Aut\ (S, o_\rho)= C_{Aut\ (S, o)} (\rho)$.
\end{prop}
\proof 
From above example it is clear that  $(S,o_\rho)$ is a twisted right gyrogroup together with unique inverse $x'_\rho= x'\theta \rho$ for each $x\in S$. Let  $h\in Aut\ (S, o)$.  Then 
\begin{eqnarray}\label{tat}
[x\ o_\rho\  y]\theta h & = &  (x\theta\rho\ o\ y)\theta h\\ \nonumber
& = & x\theta(\rho h)\  o\ y\theta h \nonumber\\
& = & x\theta (\rho h\rho)\  o_\rho\ y\theta h \nonumber
\end{eqnarray}
But
 \begin{eqnarray*}
[([x]'_\rho)\theta h]'_\rho & = & [(x'\theta \rho)\theta h]'_\rho\\
&= & [x'\theta (\rho\theta h)]'\theta\rho\\
& = & x\theta (\rho h \rho)
\end{eqnarray*}
for\ $\rho, h\in Aut \ (S, o)$ and $x'' = x $.  Thus 
$$h (xo_\rho y)= [h([x]'_\rho)]'_\rho o_\rho h(y), $$
This shows that $h\in TA (S, o_\rho)$. Using equation (\ref{tat}),  we observe that $h\in Aut\ (S, o_\rho)\Leftrightarrow h= \rho h \rho$ or equivalently $h\in C_{Aut\ (S, o)} (\rho)$.\endproof

\begin{exl}\label{4ag}
Let $S = \{e,\ x_1,\  x_2,\ \ldots, \ x_{n-1}\}$. Define a binary operation $o$ on $S$ by taking $e$ as the identity and defining $x_i \ o \ x_j = x_i$ if $i\ne j$ and $x_i\ o \ x_i = e$. Then $(S, o)$ is a right loop with $x' = x$ and $f(x', x) = I_S$ for all $x\in S$. Also, for $i\ne j$
\begin{displaymath}
x_k\theta f(x_i, x_j)  =  \left\{ 
\begin{array}{lcr}
x_k &{\rm for} & i, j\neq k\\
x_j & {\rm for} & k=i\\
x_i & {\rm for} & k=j
\end{array}
\right.
\end{displaymath}
This shows that the group torsion $G_S = Sym\ (S\setminus \{e\})$. It is also evident that $f(x_i, x_j)\in Aut\ (S, o)$. Thus, $Aut\ (S, o) = G_S$ and so $(S, o)$ is a right gyrogroup~\cite{ylal}. Take $\rho = f(x_1, x_2)$, then from Proposition~\ref{tat1}, $TAut\ (S, o_\rho) = Sym\ (S\setminus \{e\})$. Again,  by Proposition~\ref{tat1} $Aut\ (S, o_\rho)= C_{Sym\ (S\setminus \{e\})} (\rho)= Sym\ (S\setminus \{e, x_1, x_2\})$.
\end{exl}

\section{Twisted Subgroups and Twisted Right Gyrogroups}

\begin{defn}\label{tsgp} (Twisted Subgroups \cite{aschbacher})
A subset $S$ of a group $G$ is called a twisted subgroup of $G$ if it satisfies the following two conditions:\\
(i)\ $1\in S$, $1$ being the identity element of the group $G$.\\
(ii)\  If $x, y\in S$, then $xyx\in S$.
\end{defn}

\begin{defn}\label{tgyrogp} (Twisted Gyrogroups)
A twisted right gyrogroup $(S, o)$ is said to be a twisted gyrogroup if the right inner mapping  $f:S\times S\rightarrow Sym\ S$ satisfies the following condition:
\begin{eqnarray}\label{rlproperty}
f(x, y ) & = & f(x, xoy)\qquad \qquad {\rm (Right\,  loop\,  property)}
\end{eqnarray}
 for all $x, y\in S$.
\end{defn}

\begin{prop}
If $(S, o)$ is a twisted gyrogroup, then we have the following identities
\begin{eqnarray}
f(x, y)^{-1} & = & f(xoy, y')\label{tggp1}\\
f(x, y)^{-1} & = & f(xoy, x)\label{tggp2}\\
f(x, y)^{-1} & = & \eta\left( f(y', x')\right)\label{tggp3}\\
f(x, y) & = & \eta\left(f(y'ox', y')\right)\label{tggp4}\\
f(x, y)^{-1} & = & f(y, x)\label{tggp5}
\end{eqnarray}
\end{prop}
\proof
If $(S, o)$ is a twisted gyrogroup, then $S$ is a twisted gyrotransversal to the subgroup $G_S$ in the group $G_S S$. 
Taking $ y=z'$ in equation~(\ref{pp1}) and using   definition of twisted right gyrogroup, we have $f(x, z') f(xo z',z) = I_S$. Replacing $z'$ by $y$, we get $f(x,  y) f(xoy, y') = I_S$.  This gives equation~(\ref{tggp1}). Using right loop property in equation~(\ref{tggp1}), we have
\begin{eqnarray*}
f(x, y)^{-1} & = & f(xoy, y')\\
& = & f(xoy, (xoy)oy')\qquad ({\rm by \ right\  loop \ property})\\
& = & f(xoy, x)\qquad \qquad ({\rm because\ (xoy)oy'= xo(yoy')= x})
\end{eqnarray*}
This proves identity (\ref{tggp2}).

Next,  in the group $G_S S$, for $x, y\in S$,
\begin{eqnarray*}
(x y)^{-1} & = & [f(x, y) (xoy)]^{-1}\\
& = & (xoy)^{-1} {f(x, y)}^{-1}\\
& = & \eta( f(x, y)^{-1}) (xoy)^{-1}\theta {f(x, y)}^{-1}
\end{eqnarray*}
Also 
\begin{eqnarray*}
y^{-1} x^{-1} & = & y'  x'\\
& = & f(y', x') (y' o x') 
\end{eqnarray*}
Hence $\eta( f(x, y)^{-1}) = f(y', x')$. Since $\eta^2 = I_S$, we have the identity~(\ref{tggp3}).

Using identity~(\ref{tggp3}), we have 
\begin{eqnarray*}
f(x, y) & = & \eta [f(y', x')^{-1}]\\
& = & \eta [f(y' o x', x'')]\ \qquad \qquad {\rm (by\ identity~(\ref{tggp1}))}\\ 
& = & \eta [f(y' o x', x)]\\
& = & \eta [f(y' o x', (y'ox')o x)]\qquad \qquad ({\rm by\ Right \ loop\ property})\\
 & = & \eta [f(y' o x', y')]\qquad \qquad {\rm as \  (y'ox')o x= y'}
\end{eqnarray*}
This proves the identity~(\ref{tggp4}).

Next, using identity~(\ref{tggp3}), we have
\begin{eqnarray*}
f(x, y)^{-1} & = & \eta [f(y', x')]\\
& = &   \eta [f(y', y'o x')]\qquad \qquad ({\rm by\ Right \ loop\ property})\\
& = & f((y'o x')', y'')^{-1} \qquad \qquad({\rm by\ identity~(\ref{tggp3})})\\
& = & f((y'o x')', y)^{-1}
\end{eqnarray*}
Thus,  $f((y'o x')', y) =  f(x, y)$. Now, using identity (\ref{tggp4}), we have 
$$f((y'o x')', y) =  \eta [f(y' o x', y')].$$
Let $z\in S$. We have unique $y\in S$ such that $y' o x' = z'$ and so  the above equation reduces to $f(z, y) =  \eta [f(z', y')] = f(y'', z'')^{-1} = f(y, z)^{-1}$ (by identity~(\ref{tggp4})). This proves the identity~(\ref{tggp5}). 
\endproof

\begin{thm}\label{twsgpgyro}
A twisted right gyrogroup $(S, o)$ is a twisted gyrogroup if and only if $S$ is a twisted subgroup of group $G_S S$ with $f(x, y)^{-1} = f(y, x),\ \forall x, y\in S$.
\end{thm}
\proof
Let $(S, o)$ be a twisted  gyrogroup.  By Theorem~\ref{mainth}, $S$ is a twisted gyrotransversal to $G_S$ in $G_S S$.  By identity~(\ref{tggp5}), $f(x, y)^{-1} = f(y, x)$. By the right loop property  $f(x, y)= f(x, xoy)= f(xoy, x)^{-1}$. Thus, $f(x, y) f(xoy, x) = 1$. Now
\begin{eqnarray*}
x y x& = & (x y) x \\
& = & f(x,  y) [(xoy) x]\\
& = & f(x, y) f(xoy, x) ((xoy)ox)\\
& = & (xoy)o x\in S
\end{eqnarray*}
This shows that $S$ is a twisted subgroup of the group $G_S S$.

Conversely, suppose that $S$ is a twisted subgroup of $G_S S$ with $f(x, y)^{-1}= f(y, x)$,\  $x, y\in S$. Let $x, y\in S$. Then
\begin{eqnarray*}
x y x& =  & f(x, y) f(xoy, x) ((xoy)ox)
\end{eqnarray*}
Since $xyx, (xoy)ox\in S$, it follows that $f(x, y) f(xoy, x) = 1$.  Thus, 
\begin{eqnarray*}
f(x, y)& = & f(xoy, x)^{-1}\\
& = & f(x, xoy)
\end{eqnarray*}
 This shows that $(S, o) $ is a twisted gyrogroup.
\endproof

\end{document}